\documentclass[a4paper]{article}
\newcommand{\be}{\begin{equation}}
\newcommand{\ee}{\end{equation}}

\newcommand{\dz}{\wedge}

\newcommand{\ba}{\begin{array}}
\newcommand{\ea}{\end{array}}
\newcommand{\beq}{\begin{eqnarray}}
\newcommand{\eeq}{\end{eqnarray}}
\newtheorem{lm}{Lemma}
\newtheorem{thee}{Theorem}
\newtheorem{proo}{Proposition}
\newtheorem{co}{Corollary}
\newtheorem{rem}{Remark}
\newtheorem{deff}{Definition}
\newcommand{\bd}{\begin{deff}}
\newcommand{\ed}{\end{deff}}

\newcommand{\bl}{\begin{lm}}
\newcommand{\el}{\end{lm}}
\newcommand{\bp}{\begin{proo}}
\newcommand{\ep}{\end{proo}}
\newcommand{\bt}{\begin{thee}}
\newcommand{\et}{\end{thee}}
\newcommand{\bc}{\begin{co}}
\newcommand{\ec}{\end{co}}
\newcommand{\brm}{\begin{rem}}
\newcommand{\erm}{\end{rem}}

\newcommand{\der}{\mbox{d}}
\begin{document}

\thispagestyle{empty}

\title {LOCALLY SASAKIAN MANIFOLDS
\footnote{Research supported 
by 
Komitet Bada\'{n} Naukowych (Grant nr 2 P03B 017 12)}\\ 
\vskip 1.truecm
{\small {\sc Micha\l ~Godli\'nski, Wojciech Kopczy\'nski\footnote{{\small{\small {\tt e-mail: kopcz@fuw.edu.pl
}}}} , Pawe\l ~Nurowski\footnote{{\small{\small {\tt e-mail: nurowski@fuw.edu.pl}}}}}} 
\vskip -0.3truecm
{\small {\it Instytut Fizyki Teoretycznej}}\\
\vskip -0.3truecm
{\small {\it Uniwersytet Warszawski, ul. Ho\.za 69, Warszawa, Poland}}\\
}
\author{\mbox{}}
\maketitle
\begin{abstract}
We show that every Sasakian manifold in dimension $2k+1$ is locally generated 
by a free real function of $2k$ 
variables. This function is a Sasakian analogue 
of the K\"ahler potential for K\"ahler geometry. It is also shown that every 
locally Sasakian-Einstein manifold in $2k+1$ 
dimensions 
is generated by a locally K\"ahler-Einstein manifold in dimension
$2k$. 
\end{abstract}

\newpage
\noindent

\rm
\section{Introduction}
The Sasakian structure, wich is defined on an odd
dimensional manifold is, in a sense, the closest possible analog of
the K\"ahler geometry of even dimension. It was introduced by S. Sasaki
\cite{Sasaki} in 1960,  who considered it as a special kind of
contact geometries. Sasakian structure consists, in paraticular, of
the contact 1-form $\eta$ and the Riemannian metric $g$. The
differential of $\eta$ defines a 2-form, which constitutes
an analog of the fundamental form of K\"ahler geometry. \\
Sasakian geometry was primarily studied as a substructure within
the catagory of contact structures. A review of this approach can be
found in \cite{Blair,Yano}. In this letter we exploit the analogy between
Sasakian and K\"ahler geometry. We show that a well known fact that a
K\"ahler geometry can be locally generated by a K\"ahler potential has
its Sasakian counterpart. This result may be of some use in
constructing a vast family of examples of Sasakian and Sasakian-Einstein
structures.\\ 
The Sasakian and Sasakian-Einstein structures appear in physics in the
context of string theory. It turns out that a metric cone $(C({\cal S})={\bf
R}_+\times {\cal S}, \bar{g}={\rm d}r^2+r^2 g)$ over a Sasakian-Einstein 
manifold $({\cal S},g)$ is K\"ahler and Ricci flat, i.e. it constitutes a
Calabi-Yau manifold. Moreover, the Sasaki-Einstein manifolds in
dimensions $2k+1$ and Sasakian manifolds with three Sasakian structures in
dimesnion $4k+3$ are related to the Maldacena conjecture
\cite{Boga,Figueroa,Maldacena,Witten}. It turns out that they are one of
very few structures which can serve as a
compact factor $\cal S$ in (anti-de-Sitter)$\times \cal S$ background
for classical field theories which, via
the Maldacena conjecture, correspond to the
large $N$ limit of certain quantum conformal field theories.\\
A formal definition of a Sasaki manifold is as follows.\\
{\bf Definition 1}\\
Let $\cal S$ be a $(2k+1)$-dimensional manifold equipped with a structure
$(\phi,\xi,\eta,g)$ such that:\\ 
(i) $\phi$ is a (1,1) tensor field,\\
(ii) $\xi$ is a vector field,\\ 
(iii) $\eta$ is a field of an 1-form,\\ 
(iv) $g$ is a Riemannian metric.\\ Assume, in addition, that for any
vector fields $X$ and $Y$ on $\cal S$,  
$(\phi,\xi,\eta,g)$ satisfy the following algebraic conditions:
\vspace{-0.3truecm}
\begin{itemize}
\item[(1)] $\phi^2 X=-X+\eta(X)\xi$,\vspace{-0.3truecm} 
\item[(2)] $\eta(\xi)=1$,\vspace{-0.3truecm}
\item[(3)] 
$g(\phi X,\phi Y)=g(X,Y)-\eta(X)\eta(Y)$,\vspace{-0.3truecm}
\item[(4)]
$g(\xi,X)=\eta(X)$,\vspace{-0.3truecm}
\end{itemize} 
and the folllowing differential conditions:\vspace{-0.3truecm}
\begin{itemize}
\item[(5)] 
$N_\phi+{\rm d}\eta\otimes\xi=0$,\\ 
where $N_{\phi}(X,Y)=[\phi X,\phi Y]+\phi^2[X,Y]-\phi[\phi X,Y]-\phi[X,\phi
Y]$ is the Nijenhuis tensor for $\phi$,\vspace{-0.3truecm}
\item[(6)]
${\rm d}\eta(X,Y)=g(\phi X,Y)$.
\end{itemize}\vspace{-0.3truecm}
Then $\cal S$ is called a Sasakian manifold.\\
{\bf Example 1}\\
A standard example of a Sasaki manifold is an odd dimensional sphere
$${\bf S}^{2k+1}=\{~{\bf
C}^{k+1}\ni (z^1,...,z^{k+1}) ~:~
|z^1|^2+...+|z^{k+1}|^2=1~\}\subset{\bf C}^{k+1},$$ viewed as a
submanifold of ${\bf C}^{k+1}$. Let $J$ be the
standard complex structure on ${\bf C}^{k+1}$, $\tilde{g}$ the
standard flat metric on ${\bf C}^{k+1}\equiv {\bf R}^{2k+2}$, and $n$ be the unit normal
to the sphere. The vector field $\xi$ on ${\bf S}^{2k+1}$ is defined
by $\xi=-Jn$. If $X$ is tangent vector to the sphere then $JX$
uniquely decomposes onto the part paralell to $n$ and the part tangent
to the sphere. Denote this decomposition by $JX=\eta(X)n+\phi X$. This
defines 1-form $\eta$ and tensor field $\phi$ on ${\bf
S}^{2k+1}$. Denoting the restriction of $\tilde{g}$ to ${\bf
S}^{2k+1}$ by $g$ we obtain $(\phi,\xi,\eta,g)$ structure on ${\bf
S}^{2k+1}$. It is a matter of checking that this structure equips
${\bf S}^{2k+1}$ with a structure of a Sasakian-Einstein
manifold. This construction is, in a certain sense, a Sasakian
counterpart of the Fubini-Study K\"ahler structure on ${\bf CP}^k$.\\
{\bf Notation}\\
We adapt the following notation:\\
$K_{,l}$ denotes the partial derivative of a function $K$ with respect to the coordinate $z^l$.\\
Complex conjugate of an indexed quantity, e.g. $a^i_{~j}$, is usually 
denoted by a bar over it, i.e. $\overline{a^i_{~j}}$. Our notation is: 
$\overline{a^i_{~j}}= \bar{a}^{\bar{i}}_{~\bar{j}}$.\\
Symmetrized tensor products of 1-forms $\eta$ and $\lambda$ is denoted by 
$\eta\lambda=\frac{1}{2}(\eta\otimes\lambda+\lambda\otimes\eta)$.\\
{\bf The main result}\\
The pouropuse of this letter is to prove the following theorem, which
locally characterizes all Sasakian and Sasakian-Einstein manifolds.\\
{\bf Theorem.}\\
Let $\cal U$ be an open set of ${\bf C}^k\times {\bf R}$ and 
let $(z^1,z^2,...,z^k,x)$ be Cartesian coordinates in $\cal U$. Consider: \vspace{-0.3truecm} 
\begin{itemize}
\item[] a vector field $\xi=\partial_x$\vspace{-0.3truecm} 
\item[] a real-valued function $K$ on $\cal U$ such that $\xi(K)=0$\vspace{-0.3truecm} 
\item[] an 1-form \[\eta={\rm d}x+i\sum_{m=1}^k(K_{,m}{\rm d}z^m)-i\sum_{\bar{m}=1}^k(K_{,\bar{m}}
{\rm d}\bar{z}^{\bar{m}})\]\vspace{-0.3truecm} \vspace{-0.3truecm} 
\item[] a bilinear form 
\[g=\eta^2+2\sum_{m,\bar{k}=1}^kK_{,m\bar{k}}{\rm d}z^m{\rm d}\bar{z}^{\bar{k}}\]\vspace{-0.3truecm} \vspace{-0.3truecm} 
\item[] a tensor field \[\phi=-i\sum_{m=1}^k[(\partial_m-iK_{,m}\partial_x)\otimes{\rm d}z^m]+i\sum_{\bar{m}=1}^k
(\partial_{\bar{m}}+iK_{,\bar{m}}\partial_x)\otimes{\rm d}{\bar{z}}^{\bar{m}}].\]
\end{itemize} \vspace{-0.3truecm}  
I) If the function $K$ is choosen in such a way that the bilinear form $g$ has 
positive definite signature then $\cal U$ equipped with the structure 
$(\phi,\xi,\eta,g)$ is a Sasakian manifold. Moreover, every Sasakian manifold 
can locally be generated by such a function $K$.\\
II) The above Sasakian structure satisfies Einstein equation
$Ric(g)=\lambda g$ if and only if $\lambda=2k$ and the function $K$
satsifies \vspace{-0.1truecm}
$$-[\log\det(K,_{i{\bar{j}}})],_{m{\bar{n}}}=2(k+1)K,_{m{\bar{n}}}.$$
\section{Almost conatact versus Sasakian manifolds}
\noindent
{\bf Definition 2}\\
Consider $(2k+1)$-dimensional manifold $\cal S$ equipped with a structure
consisting of a (1,1) tensor field $\phi$, a vector field $\xi$ and a
field of an 1-form $\eta$. Assume, in addition, that for every vector field $X$ on
$\cal S$ $(\xi,\eta,\phi)$  satisfy the following algebraic conditions:\vspace{-0.2truecm}
\begin{itemize}
\item[(1)] $\phi^2 X=-X+\eta(X)\xi$,\vspace{-0.3truecm} 
\item[(2)] $\eta(\xi)=1$.
\end{itemize}\vspace{-0.3truecm}
Then $\cal S$ is called almost contact manifold. If, in addition, an
almost contact manifold $({\cal S},(\xi,\eta,\phi))$ is equipped with
Riemannian metric $g$ such that for every veector fields $X$ and $Y$ on
$\cal S$ we have \vspace{-0.2truecm}
\begin{itemize}
\item[(3)] 
$g(\phi X,\phi Y)=g(X,Y)-\eta(X)\eta(Y)$,\vspace{-0.3truecm}
\item[(4)]
$g(\xi,X)=\eta(X)$,
\end{itemize}\vspace{-0.3truecm}
then the almost contact manifold is called almost contact metric manifold.\\
Note that every Sasakian manifold is an almost contact metric
manifold.\\

\noindent
Let T$^{\bf C} \cal S$ be the complexification of the tangent bundle
of an almost contact manifold $\cal S$. The almost contact structure
$(\xi,\eta,\phi)$ on $\cal S$ defines the decomposition
$${\rm T}^{\bf C}{\cal S}={\bf C}\otimes \xi\oplus N\oplus\bar{N},$$
where ${\bf C}\otimes \xi$, $N$ and
$\bar{N}$ are eigenspaces of $\phi$ with eigenvalues 0, $-i$ and $i$,
respectively.\\
We say that a vector subbundle $Z$ of T$^{\bf C}\cal S$ is involutive
if and only if $[\Gamma(Z),\Gamma(Z)]\subset \Gamma(Z)$, where
$\Gamma(Z)$ denotes the set of all sections of $Z$. \\
{\bf Lemma }\\
For an almost contact structure the condition 
$N_\phi+{\rm d}\eta\otimes\xi=0$ is satisfied if and only if the
bundle $N$ is involutive, $[\Gamma(N),\Gamma(N)]\subset \Gamma(N)$,
and $[\xi,\Gamma(N)]\subset \Gamma(N)$.\\
{\bf Proof}.\\
Let $X,Y\in\Gamma(N)$. Making use of the eigenvalue property of $\phi$
and of property (1) of Definition 2, we get the following expressions
for the Nijenhuis tensor of $\phi$:
$$N_\phi(X,Y)=-2[X,Y]+2i\phi([X,Y])+\eta([X,Y])\xi,$$
$$N_\phi(X,\bar{Y})=\eta([X,\bar{Y}])\xi,$$
$$N_\phi(X,\xi)=-[X,\xi]+i\phi([X,\xi])+\eta([X,\xi])\xi.$$
Observe that the last component of the above formulae is the action of
$-{\rm d}\eta\otimes\xi$ on $(X,Y)$, $(X,\bar{Y})$ and $(X,\xi)$,
respectively. Therefore $N_\phi+{\rm d}\eta\otimes\xi=0$ if and only
if 
$$\phi([X,Y])=-i[X,Y]$$
and
$$\phi([X,\xi])=-i[X,\xi].$$
This finishes the proof.\\
{\bf Corollary}\\
For an almost contact structure satisfying $N_\phi+{\rm
d}\eta\otimes\xi=0$ (in particular for a Sasaki structure) the bundle
${\bf C}\otimes\xi\oplus N$ is involutive. 
\section{Sasakian geometry in a null frame}
Let $({\cal S},(\xi,\eta,\phi,g))$ be a Sasakian manifold of dimension
$2k+1$. The algebraic conditions (1)-(4) of Definition 1 imply an
existence of a local basis $(\xi, m_i,\bar{m}_{\bar{i}})$,
$i,\bar{i}=1,2,...,k$, of
complex-valued vector fields on $\cal S$, with a cobasis $(\eta, \mu^i,
\bar{\mu}^{\bar{i}})$, such that 
\be 
g=\eta^2+ 2 \sum_{l=\bar{l},l=1}^k \mu^l\bar{\mu}^{\bar{l}},\label{g}
\ee
\be
\phi=i\sum_{\bar{j}=1}^k(\bar{m}_{\bar{j}}\otimes\bar{\mu}^{\bar{j}})-i\sum_{j=1}^k(m_{j}\otimes
\mu^j).\label{phi}
\ee
Since $({\cal S},(\xi,\eta,\phi,g))$ is Sasakian then its bundle
${\bf C}\otimes \xi\oplus N$ is involutive. This is equivalent to
the condition that the forms $\mu^1,\mu^2,...,\mu^k$ generate a closed
differential ideal i.e. 
\be
{\rm d}\mu^i\dz \mu^1\dz \mu^2\dz ...\dz \mu^k=0\quad\quad\quad \forall
i=1,2,...,k.\label{frob}
\ee
Condition (6) of Definition 1 of a Sasakian manifold in this basis
reads
\be
{\rm d}\eta=-2i \sum_{l=\bar{l},l=1}^k \mu^l\dz \bar{\mu}^{\bar{l}}.\label{deta}
\ee
Thus, the fact that the manifold is Sasakian necessarily implies the 
existence of a local basis
$(\xi,m_i,\bar{m}_{\bar{i}})$ with a dual basis $(\eta,
\mu^i,\bar{\mu}^{\bar{i}})$ such that (\ref{g})-(\ref{deta}) holds. The
converse is also true: if a real vector field $\xi$ on a manifold
$\cal S$ can be suplemented by $k$ complex-valued vector fields $m_i$
such that $(\xi,m_i,\bar{m}_{\bar{i}})$ and
$(\eta,\mu^i,\bar{\mu}^{\bar{i}})$ form a mutually dual basis for T$^{\bf C}\cal
S$ and T$^{*\bf C}\cal S$, respectively, satisfying
(\ref{frob})-(\ref{deta}), then the structure
$(\xi,\eta,\phi,g)$ defined by $\xi$, $\eta$ and $g$, $\phi$ of 
(\ref{g})-(\ref{phi}) is a
Sasakian manifold. This fact can be seen by 
noting that condition (\ref{frob}) is equivalent to
the existence of complex-valued functions $a_{ijk}$, $b_{i{\bar{j}}k}$
and $c_{ij}$ such that 
\be 
{\rm
d}\mu^i=\sum_{j,n=1}^k a_{ijn}\mu^j\dz \mu^k+\sum_{\bar{j},n=1}^k
b_{i{\bar{j}}n}\bar{\mu}^{\bar{j}}
\dz \mu^k+\sum_j^k c_{ij}^j\dz \eta.\label{dm}
\ee
The dual conditions to conditions
(\ref{deta})-(\ref{dm}) imply that $N$ is involutive and that $[\xi, \Gamma
(N)]\subset \Gamma(N)$. These, when compared with Lemma of
Corollary 1, imply condition (7), which is the only condition from
Definition 1 which, a'priori, was not assumed for $(\xi,\eta,\phi,g)$.
This proves the following Proposition.\\
{\bf Proposition 1}\\
(I) (Local version)\\
Let $(\xi,\eta,\phi,g)$ be a Sasakian structure on a manifold $\cal S$
of dimension $2k+1$. Then there exisists a local basis
$(\xi,m_i,\bar{m}_{\bar{i}})$, $i,\bar{i}=1,2,...k$ of T$^{\bf C}\cal
S$ with dual basis $(\eta,\mu^i,\bar{\mu}^{\bar{i}})$ such that 
$$ 
g=\eta^2+ 2 \sum_{l=\bar{l},l=1}^k \mu^l\bar{\mu}^{\bar{l}}
$$
$$
\phi=i\sum_{\bar{j}=1}^k(\bar{m}_{\bar{j}}\otimes\bar{\mu}^{\bar{j}})-i\sum_{j=1}^k(m_{j}\otimes
\mu^j),$$
$$
{\rm d}\mu^i\dz \mu^1\dz \mu^2\dz ...\dz \mu^k=0\quad\quad\quad \forall
i=1,2,...,k,
$$
$${\rm d}\eta=-2i \sum_{l=\bar{l},l=1}^k \mu^l\dz \bar{\mu}^{\bar{l}}.$$
(II) (Global version)\\
Every almost contact metric structure which satisfies condition (6) of
Definition 1 is Sasakian if and only if its canonical decomposition 
${\rm T}^{\bf C}{\cal S}={\bf C}\otimes \xi\oplus N\oplus\bar{N},$
consists of involutive ${\bf C}\otimes \xi\oplus N$ part. \\

\noindent
We close this section with a quick application of part (I) of
Proposition 1. It is well known that a vector field $\xi$ on a
Sasakian manifold $({\cal S},(\xi,\eta,\phi,g))$ is a Killing vector
field. This in particular means that the Lie derivatives ${\cal
L}_\xi$ of $g$ and $\eta$ vanish. The second of these facts is an
immediate consequence of (\ref{deta}). To calculate ${\cal L}_\xi g$
one uses (\ref{deta}) and (\ref{dm}). After some work one shows that
vanishing of ${\cal L}_\xi g$ is equivalent to
$c_{ij}+\overline{c_{ji}}=0$ where $c_{ij}$ are functions appearing in
(\ref{dm}). On the other hand, these equations are automatically
implied by application of d on both sides of equation (\ref{deta}).
\section{Analogue of the K\"ahler potential}
We pass to a construction of local coordinates on a Sasakian
manifold $({\cal S},(\xi,\eta,\phi,g))$. We assume that all the fields
defining the Sasakian structure are smooth on $\cal S$.\\ 
In a considered region of $\cal S$, we chose a local frame 
$(\xi,m_i,\bar{m}_{\bar{i}})$ of Proposition 1. Now, the fact that
$\xi$ is a Killing vector field on $\cal S$ together with 
the complex version of the Fr\"obenius theorem, assures that condition 
(\ref{frob}) is equivalent to an existence of complex-valued functions
$f^i_{~j}$ and $z^i$, $i,j=1,2,...k$ such that 
\be
\mu^i=f^i_{~j}{\rm d}z^j.
\ee  
Since the forms $(\mu^i,\bar{\mu}^{\bar{i}})$ form a part of the basis on
the considered region of $\cal S$ then we also have
\be
{\rm d}z^1\dz {\rm d}z^2\dz ...{\rm d}z^k\dz{\rm d}\bar{z}^1\dz{\rm
d}\bar{z}^2\dz ...\dz{\rm d}\bar{z}^k\neq 0.
\ee
For the basis $(\xi,\partial_{z^1},...,\partial_{z^k},
\partial_{\bar{z}^1},...,\partial_{\bar{z}^k})$ and its
dual $(\eta,\der z^1,...,\der z^k,\der\bar{z}^1,...,\der\bar{z}^k)$
the Maurer-Cartan relations for $\der(\der
z^i)=0=\der(\der\bar{z}^{\bar{i}})$ readilly show that
$[\xi,\partial/\partial
z^i]=0=[\xi,\partial/\partial\bar{z}^{\bar{i}}]$. Therefore, there
exists a real coordinate $x$ complementary to
$z^1,...,z^k,\bar{z}^1,...,\bar{z}^k$ such that 
\be
\xi=\partial_x\quad
\ee
and the form $\eta$ reads
\be
\eta={\rm d}x+ p_j{\rm d}z^j+\bar{p}_{\bar{j}}{\rm
d}\bar{z}^{\bar{j}}.
\ee
Comparing this with the fact that $\xi$ preserves $\eta$ leads to 
the conclusion that the functions $p_i$ are independent of coordinate
$x$, $\partial p_i/\partial x=0$.\\
Condition (\ref{deta}) is now equivalent to the follwing two
conditions for the differentails of functions $p_i$:
\be
p_{i,j}-p_{j,i}=0\label{pij}
\ee
and
\be
p_{j,\bar{i}}-\bar{p}_{\bar{i},j}=2i\sum_{l=\bar{l}, l=1}^k 
f^l_{~j}\bar{f}^{\bar{l}}_{~\bar{i}}.\label{vij}
\ee
In a simply connected region of $\cal S$ equation (\ref{pij})
guarantees an existence of a complex-valued function $V$ such that
\be 
p_i=\frac{\partial V}{\partial z^i}.
\ee 
Since $p_i$ is independent of $x$ it is
enough to consider functions $V$ such that $\partial V/\partial
x=0$. Inserting so determined $p_i$ to equation (\ref{vij}) we show
that now (\ref{vij}) is equivalent to 
\be
K_{,j\bar{i}}=\sum_{l=\bar{l}, l=1}^k 
f^l_{~j}\bar{f}^{\bar{l}}_{~\bar{i}}, \label{kij}
\ee
where we have introduce Im$V=K$ and Re$V=L$. Finally we note that now
$$
\eta={\rm d}(x+L)+i\sum_{j=1}^k K_{,j}{\rm d}z^j-i\sum_{\bar{j}=1}^k
K_{,\bar{j}}{\rm d}\bar{z}^{\bar{j}},
$$
so redefining the $x$ coordinate by $x\to x+L$ we simplify $\eta$ to
the form $\eta={\rm d}x+i\sum_{j=1}^k K_j{\rm d}z^j-i\sum_{\bar{j}=1}^k
K_{\bar{j}}{\rm d}\bar{z}^{\bar{j}}$.
Using (\ref{kij}) we can eliminate functions $f^i_{~j}$ from
formulae defining our Sasakian structure. Indeed, 
$$
g=\eta^2+ 2 \sum_{l=\bar{l},l=1}^k \mu^l\bar{\mu}^{\bar{l}}=$$\vspace{-0.5truecm}
$$
\eta^2+ 2 \sum_{l=\bar{l},l=1}^k\sum_{j,\bar{i}=1}^k
f^l_{~j}\bar{f}^{\bar{l}}_{~\bar{i}}{\rm d}z^j{\rm
d}\bar{z}^{\bar{i}}=$$\vspace{-0.2truecm}
$$\eta^2+ 2 \sum_{j,\bar{i}=1}^k
K_{,j\bar{i}}{\rm d}z^j{\rm
d}\bar{z}^{\bar{i}}.$$
In this
way we obtain the following theorem.\\
{\bf Theorem 1.}\\
Let $\cal U$ be an open set of ${\bf C}^k\times {\bf R}$ and 
let $(z^1,z^2,...,z^k,x)$ be Cartesian coordinates in $\cal U$. Consider: \vspace{-0.3truecm} 
\begin{itemize}
\item[(i)] a vector field $\xi=\partial_x$\vspace{-0.3truecm} 
\item[(ii)] a real-valued function $K$ on $\cal U$ such that $\xi(K)=0$\vspace{-0.3truecm} 
\item[(iii)] an 1-form \[\eta={\rm d}x+i\sum_{m=1}^k(K_{,m}{\rm d}z^m)-i\sum_{\bar{m}=1}^k(K_{,\bar{m}}
{\rm d}\bar{z}^{\bar{m}})\]\vspace{-0.3truecm} \vspace{-0.3truecm} 
\item[(iv)] a bilinear form 
\[g=\eta^2+2\sum_{m,\bar{k}=1}^kK_{,m\bar{k}}{\rm d}z^m{\rm d}\bar{z}^{\bar{k}}\]\vspace{-0.3truecm} \vspace{-0.3truecm} 
\item[(v)] a tensor field \[\phi=-i\sum_{m=1}^k[(\partial_m-iK_{,m}\partial_x)\otimes{\rm d}z^m]+i\sum_{\bar{m}=1}^k
(\partial_{\bar{m}}+iK_{,\bar{m}}\partial_x)\otimes{\rm d}{\bar{z}}^{\bar{m}}].\]
\end{itemize} \vspace{-0.3truecm}  
If the function $K$ is choosen in such a way that the bilinear form $g$ has 
positive definite signature then $\cal U$ equipped with the structure 
$(\phi,\xi,\eta,g)$ is a Sasakian manifold. Moreover, every Sasakian manifold 
can locally be generated by such a function $K$.\\

\noindent
The function $K$ appearing in the above theorem is a Sasakian analogue
of the K\"ahler potential generating K\"ahler geometries. We call it a
Sasakian potential.\\

\noindent
We close this section with a remark that several Sasakian potentials
may generate the same Sasakian structure. This is evident if one notes
that the following transformations 
\be
K\to K+f(z^j)+\bar{f}(\bar{z}^{\bar{j}})\quad\quad\quad x\to x+
i\bar{f}(\bar{z}^{\bar{j}})-if(z^j),\label{trans}
\ee 
with $f$ beeing a holomorphic function of $z^j$s, do not change the
Sasakian structure of Theorem 1. Thus, transformations (\ref{trans})
are the gauge transformations for the Sasakian potential.
\section{Locally Sasakian-Einstein structures} 
In this section we calculate the Ricci tensor for the Sasakian metric 
$g$ generated in a region $\cal U$ by the Sasakian potential $K$ of 
Theorem 1. We also derive the equation which the Sasakian potential
has to obey for Sasakian metric to satisfy Einstein equations
$Ric(g)=\lambda g$. In this section we use the Einstein summation 
convention.\\
Let $\cal U$ be a simply connected region of of ${\bf C}^k\times {\bf
R}$ as in Theorem 1. Consider a Sasakian structure defined in this
Theorem by the Sasakian potential $K$. In the holonomic cobasis $$({\rm
d}y^\mu)=({\rm
d}x, {\rm d}z^j,{\rm d}\bar{z}^{\bar{j}})$$ 
the covariant components of the Sasakian metric
read
\be
g_{\mu\nu}=\left (
\begin{array}{ccc}
1&iK_{,j}&-iK_{,\bar{j}}\\
iK_{,i}&-K_{,i}K_{,j}&K_{,i\bar{j}}+K_{,i}K_{,\bar{j}}\\
-iK_{,\bar{i}}&K_{,\bar{i}j}+K_{,\bar{i}}K_{,j}&-K_{,\bar{i}}K_{,\bar{j}}
\end{array}\right ).
\ee
The contravariant components of the metric read
\be
g^{\mu\nu}=\left (
\begin{array}{ccc}
1+2K_{,i}K_{,\bar{j}}\kappa^{i\bar{j}}&iK_{,\bar{j}}\kappa^{i\bar{j}}&-iK_{,j}\kappa^{j\bar{i}}\vspace{0.2truecm}\\
iK_{,\bar{i}}\kappa^{j\bar{i}}&0&\kappa^{j\bar{i}}\\
-iK_{,i}\kappa^{i\bar{j}}&\kappa^{i\bar{j}}&0
\end{array}\right ),
\ee
where
\be
\kappa^{j\bar{l}}K_{,i\bar{l}}=\delta^j_{~i}\quad\quad\kappa^{l\bar{j}}K_{,l\bar{i}}
=\delta^{\bar{j}}_{~\bar{i}}\quad\quad\overline{\kappa^{l\bar{j}}}=\kappa^{j\bar{l}}=\kappa^{\bar{l}j}.
\ee
The connection 1-forms
$\Gamma_{\mu\nu}=\frac{1}{2}(g_{\mu\nu,\rho}+g_{\rho\mu,\nu}-g_{\nu\rho,\mu}){\rm
d}y^\rho$ read
$$
\Gamma_{xx}=0\quad\quad\quad\Gamma_{xi}=iK_{,ij}{\rm
d}z^j\quad\quad\quad\Gamma_{x\bar{i}}=\overline{\Gamma_{xi}}
\quad\quad\quad\Gamma_{ix}=iK_{,i\bar{j}}{\rm d}\bar{z}^{\bar{j}}
\quad\quad\quad\Gamma_{\bar{i}x}=\overline{\Gamma_{ix}}
$$
\be
\Gamma_{ij}=-K_{,i}K_{,jl}{\rm d}z^l-K_{,j}K_{,i\bar{l}}{\rm
d}\bar{z}^{\bar{l}}
\quad\quad\quad\Gamma_{\bar{i}\bar{j}}=\overline{\Gamma_{ij}}
\ee
$$
\Gamma_{i\bar{j}}=iK_{,i\bar{j}}{\rm d}x-K_{,l}K_{,i\bar{j}}{\rm
d}z^l+
(K_{,i\bar{j}\bar{l}}+K_{,i}K_{,\bar{j}\bar{l}}+K_{,\bar{j}}K_{,i\bar{l}}+
K_{,\bar{l}}K_{,i\bar{j}}){\rm
d}\bar{z}^{\bar{l}}\quad\quad\quad\Gamma_{,\bar{i}j}=
\overline{\Gamma_{i\bar{j}}}.
$$
It is convenient to introduce the following functions:
$$
C^i_{~jm}=\kappa^{i\bar{l}}K_{,\bar{l}jm}\quad\quad\quad B^i_{~jm}=C^i_{~jm}+
\delta^i_{~m}K_{,j}+\delta^i_{~j}K_{,m}\quad\quad\quad A_{jm}=C^l_{~jm}K_{,l}+
2K_{,j}K_{,m}-K_{,jm}
$$
$$
C^{\bar{i}}_{~\bar{j}\bar{m}}=\overline{C^i_{~jm}}\quad\quad\quad 
B^{\bar{i}}_{~\bar{j}\bar{m}}=\overline{B^i_{~jm}}\quad\quad\quad 
A_{~\bar{j}\bar{m}}=\overline{A_{~jm}}.
$$
Then the connection 1-forms $\Gamma^\mu_{~\nu}$ read
$$
\Gamma^x_{~x}=-{\rm d}K\quad\quad\quad\Gamma^x_{~j}=-K_{,j}{\rm d}x-
iA_{jm}{\rm d}z^m\quad\quad\quad\Gamma^i_{~x}=-i{\rm d}z^i\quad\quad\quad
\Gamma^x_{~\bar{j}}=\overline{\Gamma^x_{~j}}\quad\quad\quad
\Gamma^{\bar{i}}_{~x}=\overline{\Gamma^i_{~x}}
$$
$$
\Gamma^i_{~j}=-i\delta^i_{~j}{\rm d}x-\delta^i_{~j}K_{,\bar{l}}{\rm d}\bar{z}^{\bar{l}}+B^i_{~jl}{\rm d}z^l\quad\quad\quad\Gamma^i_{~\bar{j}}=-K_{,\bar{j}}{\rm d}z^i\quad\quad\quad\Gamma^{\bar{i}}_{~\bar{j}}=\overline{\Gamma^i_{~j}}
\quad\quad\quad\Gamma^{\bar{i}}_{~j}=\overline{\Gamma^i_{~\bar{j}}}.
$$
The curvature 2-forms $\Omega^\mu_{~\nu}=\frac{1}{2}R^\mu_{~\nu\rho\sigma}
{\rm d}y^\rho\dz{\rm d}y^\sigma=
{\rm d}\Gamma^\mu_{~\nu}+\Gamma^\mu_{~\rho}\dz\Gamma^\rho_{~\nu}$ read
$$
\Omega^x_{~x}=iK_{,l}{\rm d}x\dz{\rm d}z^l-iK_{,\bar{l}}{\rm d}x\dz{\rm d}\bar{z}^{\bar{l}}
$$
$$
\Omega^x_{~j}=-K_{,j}K_{,l}{\rm d}x\dz{\rm d}z^l+
(K_{,j\bar{l}}+K_{,j}K_{,\bar{l}}){\rm d}x\dz{\rm d}\bar{z}^{\bar{l}}+iA_{jl,\bar{m}}{\rm d}z^l\dz{\rm d}\bar{z}^{\bar{m}}
$$
$$
\Omega^i_{~x}=-\delta^i_{~j}{\rm d}x\dz{\rm d}z^j+i\delta^i_{~j}K_{,l}{\rm d}z^j\dz{\rm d}z^l+iK_{,\bar{l}}{\rm d}\bar{z}^{\bar{l}}\dz{\rm d}z^i
$$
$$
\Omega^i_{~j}=-i\delta^i_{~l}K_{,j}{\rm d}x\dz{\rm d}z^l+
(B^i_{~mn}B^m_{~jl}-B^i_{~jn,l}+A_{jn}\delta^i_{~l}){\rm d}z^n\dz{\rm d}z^l+
(K_{,j}K_{,\bar{l}}\delta^i_{~n}-\delta^i_{~j}K_{,n\bar{l}}-B^i_{~jn,\bar{l}})
{\rm d}z^n\dz{\rm d}\bar{z}^{\bar{l}}
$$
$$
\Omega^i_{~\bar{j}}=i\delta^i_{~l}K_{,\bar{j}}{\rm d}x\dz{\rm d}z^l+
(K_{,\bar{j}l}+K_{,\bar{j}}K_{,l})\delta^i_{~n}{\rm d}z^n\dz{\rm d}z^l-\delta^i_{~n}K_{,\bar{j}}K_{\bar{l}}{\rm d}z^n\dz{\rm d}\bar{z}^{\bar{l}}
$$
$$
\Omega^x_{~\bar{j}}=\overline{\Omega^x_{~j}}\quad
\Omega^{\bar{i}}_{~x}=\overline{\Omega^i_{~x}}\quad
\Omega^{\bar{i}}_{~\bar{j}}=\overline{\Omega^i_{~j}}\quad
\Omega^{\bar{i}}_{~j}=\overline{\Omega^i_{~\bar{j}}}.
$$
The Ricci tensor $R_{\nu\sigma}=R^\mu_{~\nu\mu\sigma}$ componenets read
$$
R_{xx}=2k\quad\quad R_{xj}=2ikK_{,j}\quad\quad R_{ij}=-2kK_{,i}K_{,j}
\quad\quad R_{\bar{i}j}=2kK_{,\bar{i}}K_{,j}-2K_{,\bar{i}j}-C^{\bar{m}}_{~\bar{m}\bar{i},j}
$$
$$
R_{x\bar{j}}=\overline{R_{xj}}\quad\quad R_{\bar{i}\bar{j}}=
\overline{R_{ij}}.
$$
Now, the Einstein equations $Ric(g)=\lambda g$, which are nontrivial
only for the components $R_{xx}$ and $R_{\bar{i}j}$ become 
$$
\lambda=2k\quad\quad\quad -(\kappa^{\bar{m}l}K_{,\bar{m}l\bar{i}})_{,j}=2(k+1)K_{,\bar{i}j}.
$$
Since the matrix $(\kappa^{\bar{m}l})$ is the inverse of
$(K_{,\bar{i}j})$ then the left hand side of the second equations
above is
$$
-(\kappa^{\bar{m}l}K_{,\bar{m}l\bar{i}})_{,j}=-\log (\det (K_{,m\bar{n}}))_{,\bar{i}j},
$$
see e.g. \cite{Landau}.\\
Thus we arrive to the following theorem.\\
{\bf Theorem 2}\\
Any Sasakian manifold of dimension $(2k+1)$ can be locally represented
by the Sasakian potential $K$ of Theorem 1. In the region where the
potential is well defined the manifold satisfies Einstein equations
$Ric(g)=\lambda g$ if and only if the cosmological constant\vspace{-0.2truecm}
$$\lambda=2k\vspace{-0.3truecm}$$
and the potential satisfies 
\be
 -\log (\det
(K_{,m\bar{n}}))_{,\bar{i}j}=2(k+1)K_{,\bar{i}j}.\label{einst}
\vspace{0.5truecm}
\ee
\noindent
Surprisingly equation (\ref{einst}) is the same as the Einstein
condition $Ric(h)=2(k+1)h$ for the K\"ahler metric 
$h=2K_{,\bar{i}j}{\rm d}\bar{z}^{\bar{i}}{\rm d}z^j$ in dimension $2k$.
Thus we
have the following Corollary.\\
{\bf Corollary}\\
Every Sasakian-Einstein manifold in dimension $(2k+1)$ is locally in
one to one correspondence with a K\"ahler-Einstein manifold in
dimension $2k$ whose cosmological constant $\lambda=2(k+1)$. The
correspondence is obtained by identifying the K\"ahler potential for
the K\"ahler-Einstein manifold with the Sasaki potential for the
Sasaki-Einstein manifold. \\

\noindent
{\bf Examples}\\
1). Sasakian potential for the sphere ${\bf S}^{2k+1}$.\\
Consider a function 
$${\cal K}=\frac{1}{2}\log (z^1\bar{z}^1+...+z^{k+1}\bar{z}^{k+1})$$
defined on ${\bf C}^{k+1}-\{0\}$. 
Let 
$$N=i({\cal K}_{,j}\der z^j-{\cal K}_{,\bar{j}}\der\bar{z}^{\bar{j}}),$$
$$H=2{\cal K}_{,i\bar{j}}\der z^i\der \bar{z}^{\bar{j}}$$
and
$$G=N^2+H.$$
The tensor fields $N$ and $G$ restrict to a sphere 
$${\bf S}^{2k+1}=\{ (z^1,...,z^{k+1})\in {\bf C}^{k+1}-\{0\}~ |~
z^1\bar{z}^1+...+z^{k+1}\bar{z}^{k+1}=1\}.$$ Denote these restrictions by 
$\eta$ and $g$, respectively. Then the 1-form $\eta$ and the Riemannian 
metric $g$ define a Sasakian-Einstein structure on ${\bf S}^{2k+1}$. 
This structure coincides with the one defined in Example 1 of Section 1. 

\noindent
To see this, recal the Hopf fibration   
${\bf U}(1)\to{\bf S}^{2k+1}\to{\bf CP}^k$ with the action of 
${\rm e}^{i\phi}\in{\bf U}(1)$ 
on $(z^1,...,z^{k+1})\in{\bf S}^{2k+1}$ defined by 
${\rm e}^{i\phi}(z^1,...,z^{k+1})=({\rm e}^{i\phi}z^1,...,{\rm e}^{i\phi}z^{k+1})$. The canonical projection is given by  
${\bf S}^{2k+1}\ni(z^1,...,z^{k+1})\to 
{\rm dir}(z^1,...,z^{k+1})\in{\bf CP}^k $. The sphere is covered by 
$k+1$ charts $$U_j=V_j\times {\bf U}(1),$$ where 
$$V_j=\{{\rm dir}(z^1,...,z^{k+1})~|~
(z^1,...,z^{k+1})\in {\bf S}^{2k+1}~{\rm and}~z^j\neq 0\}.$$ The local 
coordinates on each $U_j$ are 
$$(\xi^{ij}=\frac{z^i}{z^j},~\phi_j=\frac{i}{2}\log\frac{z^j}{\bar{z}^j}),~~ 
i=1,...,k+1, ~i\neq j.$$
Then on each chart $U_j$ the form $\eta_j=\eta|_{U_j}$ reads
$$\eta_j=\der\phi_j+\frac{i}{2}\frac{\sum_{i=1,i\neq j}^{k+1}(\bar{\xi}^{ij}\der\xi^{ij}-\xi^{ij}\der
\bar{\xi}^{ij})}{1+\sum_{i=1,i\neq j}^{k+1}|\xi^{ij}|^2}.$$
The metric $g$ restricted to $U_j$ is 
$$g_j=(\eta_j)^2+\frac{(1+\sum_{i=1,i\neq j}^{k+1}|\xi^{ij}|^2)
(\sum_{i=1,i\neq j}^{k+1}|\der\xi^{ij}|^2)-|\sum_{i=1,i\neq
j}^{k+1}(\xi^{ij}\der\bar{\xi}^{ij})|^2}{(1+\sum_{i=1,i\neq
j}^{k+1}|\xi^{ij}|^2)^2}.$$
Now, observe that on each $U_j$ the structure 
$(g_j,\eta_j)$ may be obtained by means of Theorems 1 and 2 choosing a
Sasakian potential 
$$K^j=\frac{1}{2}\log(1+\sum_{i=1,i\neq
j}^{k+1}|\xi^{ij}|^2)$$
on the corresponding $V_j$. It is easy to check that $K^j$
satisfies equation (\ref{einst}) on $V_j$. Thus, Theorem 2 assures
that the Sasakian structure generated by $(g_j,\eta_j)$ is 
Sasakian-Einstein. Easy,
but lenghty, calculation shows that the Weyl tensor of $g_j$ vanishes
identically on $U_j$. This proves that $(U_j,g_j)$ is locally 
isometric to a standard Riemannian structure on ${\bf
S}^{2k+1}$. Since $(g_j,\eta_j)$ originate from the global structure 
$(g,\eta)$ then this global Sasakian ~structure must ~coincide with the~
standard Sasakian ~structure of Example 1. Note also that $h_j=g_j-(\eta_j)^2$
projects to $V_j$ and patched together defines the Fubini-Study metric
on ${\bf CP}^k$. 
In this sense the standard Sasakian structure on ${\bf S}^{2k+1}$ 
described in Example 1 is the analogue of the 
Fubini-Study K\"ahler structure on ${\bf CP}^k$.\\  
2.) Sasakian-Einstein structure on ${\bf C}^q\times {\bf C}^n\times\bf R$.\\
Consider a function 
$$K=\frac{1}{q+n+1}[\sum_{i=1}^q\log(1+|v^i|^2)]+\frac{n+1}{2(q+n+1)}
\log(1+\sum_{I=1}^n|w^I|^2)$$ 
defined on ${\bf C}^q\times{\bf C}^n$, with coordinates 
$(z^\mu)=(v^i,w^I)$. It is easy to check that 
$$[\log\det
(K_{,\mu\bar{\nu}})]_{,\rho\bar{\sigma}}=-2(q+n+1)K_{,\rho\bar{\sigma}}.$$
Thus, via Theorems 1 and 2, such $K$ generates a Sasakian-Einstein
structure on ${\bf C}^q\times{\bf C}^n\times {\bf R}^1.$ \\
3) Locally Sasakian-Einstein structures in dimension 5.\\
If $k=2$ then, modulo the gauge (\ref{trans}), equation (\ref{einst})
may be integrated to the form
$$K_{,1\bar{1}}K_{,2\bar{2}}-K_{,1\bar{2}}K_{,2\bar{1}}={\rm e}^{-6K}.$$
This is a well known equation describing the gravitational instantons
in four dimensions \cite{BFP,Gibbons,Przan,Przan2}. Examples of the
K\"ahler-Einstein metrics generated by its solutions can be found
e.g. in  \cite{Gibbons,Pedersen,Plebprzan}. Via Theorems 1 and 2, each of these
K\"ahler-Einstein structures defines a nontrivial 
Sasakian-Einstein manifold in dimension 5.

\end{document}